\documentclass[a4paper,12pt,reqno]{article}

\usepackage[hideboxes]{boxedeps}

\usepackage{latexsym}

 \usepackage{amsmath}
\usepackage{hhline}

\newcommand{\pic}[2]{\BoxedEPSF{#1 scaled #2}}

\newcommand{\ints}{\mathbf{Z}}

\newcommand{\lm}{(\lambda,\mu)}
\newcommand{\la}{\lambda}
\newcommand{\La}{\Lambda}

\newcommand{\de}{\delta}
\newcommand{\De}{\Delta}
\newcommand{\si}{\sigma}
\newcommand{\be}{\beta}
\newcommand{\ph}{\varphi}
\newcommand{\qlm}{Q_{\lambda,\mu}}
\newcommand{\CS}{\mathcal{S}}
\newcommand{\CC}{\mathcal{C}}
\newcommand{\CA}{\mathcal{A}}
\newcommand{\phia}{\varphi}
\newcommand{\phib}{\bar{\varphi}}
\newcommand{\del}{\frac{v^{-1}-v}{s-s^{-1}}}
\newcommand{\mirror}{\overline{\phantom{oo}}}
\newcommand{\beqn}{\begin{eqnarray*}}
\newcommand{\eeqn}{\end{eqnarray*}}
\newcommand{\beqnn}{\begin{eqnarray}}
\newcommand{\eeqnn}{\end{eqnarray}}

\newcommand{\x}{\times}

\def\Xor{\pic{xor.ART} {400}}
\def\Yor{\pic{yor.ART} {400}}
\def\Ior{\pic{ior.ART} {400}}
\def\Sigmaior{\pic{sigmaior.ART} {400}}
\def\Loop{\pic{loop.ART} {400}}

\def\Idor{\pic{idor.ART} {400}}
\def\Rcurlor{\pic{rcurlor.ART} {400}}

\def\Am{\pic{am.ART} {600}}

\def\Gnplus{\pic{gnplus.ART} {600}}
\def\Gshift{\pic{gn.ART} {600}}
\def\Gcross{\pic{gplus.ART} {400}}
\def\Yn{\pic{yn.ART} {700}}

\def\Aid{\pic{ea.ART} {500}}
\def\Aa{\pic{a.ART} {500}}
\def\Aaa{\pic{aa.ART} {500}}
\def\Aainv{\pic{ainv.ART} {500}}
\def\Annulus{\pic{fa.ART} {500}}
\def\LX{\pic{lx.ART} {500}}
\def\RX{\pic{rx.ART} {500}}

\def\XC{\pic{xc.ART} {600}}
\def\XCYC{\pic{xcyc.ART} {600}}
\def\XYC{\pic{xyc.ART} {600}}

\def\Fnp{\pic{fnp.ART} {400}}
\def\TC{\pic{closeT.ART} {600}}

\def\fourtwotwo{\pic{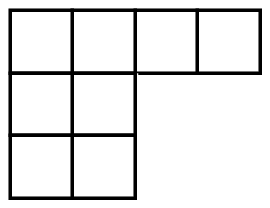} {500}}
\def\threetwo{\pic{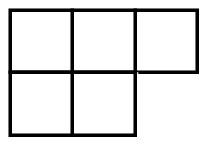} {500}}

\def\Xone{\pic{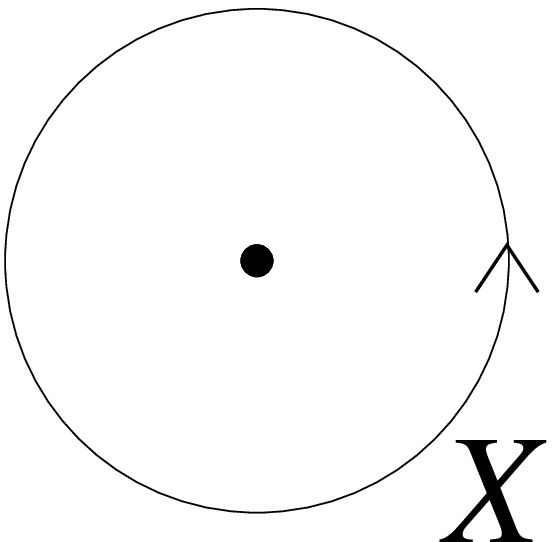} {200}}
\def\Xtwo{\pic{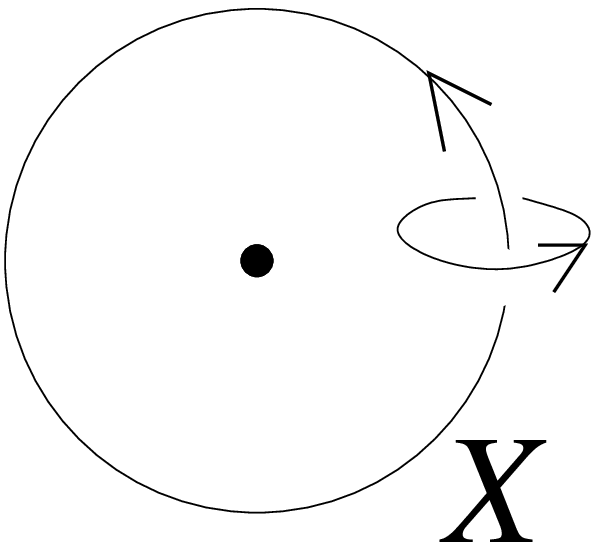} {200}}
\def\Xthree{\pic{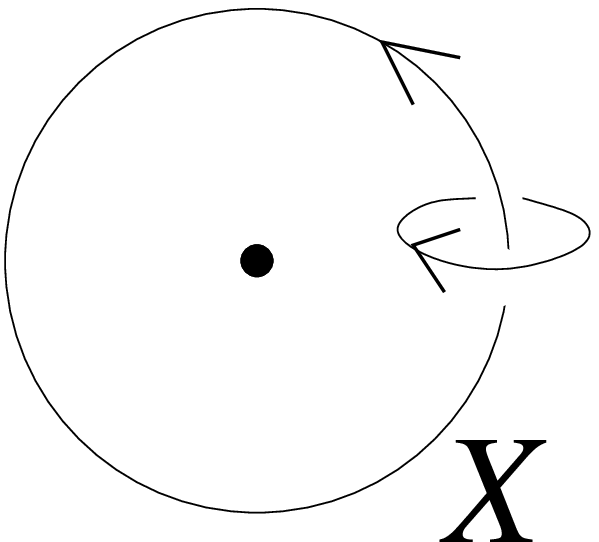} {200}}

\newtheorem{theorem}{Theorem}[section]
\newtheorem{corollary}[theorem]{Corollary}
\newtheorem{lemma}[theorem]{Lemma}

\newenvironment{remark}{\par\smallskip%
\noindent\textbf{Remark.}\  }%
{\par\smallskip}

\newenvironment{example}{\par\smallskip%
\noindent\textbf{Example.}\  }%
{\par\smallskip}

\newenvironment{Def}{\par\smallskip%
\noindent\textbf{Definition.}\  }%
{\par\smallskip}

\newenvironment{definition}{\par\smallskip%
\noindent\textbf{Definition.}\  }%
{\par\smallskip}

\newenvironment{notation}{\par\smallskip%
\noindent\textbf{Notation.}\  }%
{\par\smallskip}

\newenvironment{proof}[1][{Proof}]{\par\smallskip%
\noindent\textit{#1: }\  }
{\hfill$\Box$\par\smallskip}

\begin{document}
\begin{center}
      {\Large\bf A basis for the full Homfly skein of the annulus}

\bigskip
RICHARD J. HADJI 
\footnote{The first author was supported by EPSRC grant 99801479}
and HUGH R. MORTON 

\medskip
 {\it Department of Mathematical Sciences, University of Liverpool,\\
     Peach St, Liverpool, L69 7ZL, England. }\\
{\tt rjhadji@o2.co.uk, morton@liv.ac.uk}
\end{center}
\vskip 1cm

\begin{abstract}
 A basis, denoted $\{\qlm\}$, for the
full Homfly skein of the annulus $\CC$ was introduced in \cite{MH02},  where
$\la$ and
$\mu$ are partitions of integers $n$ and $p$ into $k$ and $k^*$ parts
respectively.  The basis consists of eigenvectors of the two meridian
maps on
$\CC$; these maps are the  linear endomorphisms of $\CC$  induced by the
insertion of a meridian loop with either orientation around a diagram in
the annulus.

In this paper we give an explicit expression for each $\qlm$ as the
determinant of a
 $(k^*+k)\x(k^*+k)$ matrix whose entries are simple 
elements $h_n, h_n^*$ in the skein $\CC$. In the
case $p=0$ ($\mu=\phi$) the determinant gives the Jacobi-Trudy formula for
the Schur function
$s_\la$ of $N$ variables as a polynomial in the complete symmetric
functions $h_n$ of the variables, \cite{Macdonald}. The Jacobi-Trudy
determinants have previously been used  by Kawagoe
\cite{Kawagoe} and Lukac \cite{Lukac} in discussing the elements in the
skein of the annulus represented by closed braids in which all strings are
oriented in the same direction. Our results and techniques here form a
natural extension of the work of Lukac.

\end{abstract}

\section*{Introduction}

The Homfly polynomials of satellites of a knot $K$ have provided an
extensive source of invariants of $K$ since the discovery of the Homfly
polynomial in 1984. It has proved helpful to organise these invariants by
means of the Homfly skein of the annulus, $\CC$, which consists of linear
combinations of diagrams in the annulus, modulo the  Homfly skein
relations, and encapsulates the relevant information about the decorating
curves used in constructing a satellite.

The skein $\CC$, which carries a natural product making it a commutative
algebra, has been studied widely following the initial work of Turaev
\cite{Turaev}, who showed it to be a polynomial algebra with an explicit
 set of  generators consisting of simple closed braids oriented in either
direction around the annulus. Other  related invariants of $K$
emerged subsequently, using the quantum groups $sl(N)_q$. The relation
between the invariants is particularly close on the
subalgebra
$\CC^+
\subset
\CC$ generated by  closed braids all oriented in the same direction. In
this context Wenzl
\cite{Wenzl} and Murakami \cite{Murakami} gave details of  elements
of
$\CC^+$ derived from the closure of idempotents of Hecke algebras, each
determined by a partition $\la$, which provide a direct translation between
Homfly invariants and  the  $sl(N)_q$ invariants arising from the
irreducible representation given by $\la$.  More explicitly skein-based
versions of the idempotents and the corresponding elements $Q_\la\in\CC^+$
were developed by Aiston
\cite{Aiston, AistonMorton}.

It has more recently been  observed, \cite{Kawagoe, Lukac}, that the
elements $Q_\la$ can be naturally characterised as eigenvectors of the
meridian maps of $\CC$, which are linear endomorphisms of $\CC$ induced at
the level of diagrams in the annulus by inserting a meridian curve around a
given diagram. The work in  \cite{Kawagoe, Lukac} shows how to describe
$Q_\la$ as the  determinant of a matrix with entries drawn from a
sequence of elements $h_n$ in $\CC$. In interpretations of $\CC^+$ as
 the representation ring of the quantum groups $sl(N)_q$ the elements
$h_n$ correspond to the irreducible representations indexed by partitions
with a single part $n$, and the determinants used give the irreducible
representation corresponding to the partition $\la$ in terms of the
Jacobi-Trudy formula, \cite[I(3.4)]{Macdonald}. 

In the full skein $\CC$ the eigenspaces of the meridian maps are all known
 to be 1-dimensional \cite{MH02}. We present here in section \ref{sect3} an
explicit determinantal formula for the corresponding eigenvectors $\qlm$,
using as matrix entries the elements
$h_n$ and the corresponding elements $h_n^*$ where the string orientation
is reversed. As in \cite{Lukac} we use an auxiliary skein $\CA$,
also a commutative algebra, based on diagrams in the annulus
which include one arc joining a point on each boundary component.

In section 1 we introduce the basic features of Homfly skein theory
needed, including properties of the skeins $\CC$ and $\CA$ and the
elements $h_n$. In section 2 we describe the sort of matrices with entries
in
$\CA$ or $\CC$ which are used, and deduce some properties of their
determinants from simple skein relations in section 1. In the final
section we construct the basis elements $\qlm$ and prove that they are
eigenvectors of the meridian maps. We show how they appear naturally as
eigenvectors when considering the Homfly satellite invariants of a framed
knot $K$, with eigenvalues that give $(1,1)$-tangle invariants of $K$. We
also oberve that the structure constants in the algebra
$\CC$ with basis $\qlm$ are non-negative integers.

Much of this work has been developed from a part of the first author's PhD
thesis
\cite{Hadji2}.

\section{The skein models}
The account here largely follows those of \cite{AistonMorton}, \cite{Lukac}
and \cite{Murphy}.

A {\em tangle} in  a planar
surface $F$, with some designated input and
output boundary points, consists of
oriented arcs in $F$ joining input points to output points and possibly some
further oriented closed curves, up to Reidemeister moves II and III. The
term {\em diagram} is often used for a tangle when there are no input or
output points.
 The {\em framed Homfly skein} ${\CS}(F)$ of $F$ is defined to be
$\Lambda$-linear combinations of oriented tangles in $F$, modulo
the two local relations 
\begin{eqnarray*}
\Xor\  -\ \Yor& =&(s-s^{-1})\ \Ior\\
\mbox{and }\  \Rcurlor &=&v^{-1}\ \Idor.
\end{eqnarray*}
 The {\em coefficient ring} can be taken as $\La=\ints[v^{\pm 1},s^{\pm 1}]$
with the elements $\{k\}=s^k-s^{-k}$ admitted as denominators for $k\ge1$.

The local relation 
\beqn \Loop\ &=&\de\ \Idor
\eeqn
is
a consequence of the main relations.
It allows the
removal of a null-homotopic closed curve without crossings, at the expense
of multiplication by the scalar $\de=\del$, except in the case where
removal of the curve leaves the empty diagram. The result can be extended to
this case too, without inconsistency, provided that the empty diagram is 
admitted  when $F$ has no designated boundary points. 

\subsection{The plane}
When $F={\bf R}^2$ every element can be represented uniquely as a scalar
multiple of the empty diagram. For a diagram $D$ the resulting scalar
$<D>\in \Lambda$ is the {\em framed Homfly polynomial} of $D$. The
more traditional Homfly polynomial $P$ is defined as the {\em ambient}
isotopy invariant which satisfies the local relation \[v^{-1}\Xor\ -\
v\Yor=(s-s^{-1})\ \Ior,\] and which takes the value $1$ on the
unknot. The framed Homfly polynomial of $D$ satisfies
\[<D>=v^{-wr(D)}\de P(D),\]
 where $wr(D)$ is the writhe of the diagram $D$.

\subsection{The Hecke algebras and extended variants}

Write $H_n$ for the skein ${\cal S}(F)$ of
$n$-tangles,  where $F$
is a rectangle with $n$ inputs at the bottom and
$n$ outputs at the top. Composing $n$-tangles by placing one above another
induces a product which makes $H_n$ into
an algebra. It has a linear basis of $n!$
elements, and is isomorphic to the Hecke algebra
$H_n(z)$, with coefficients extended to the ring $\Lambda$. This algebra has a
presentation generated by the  elementary
braids
\[\sigma_i\ =\ \Sigmaior\] subject to 
 the
braid relations
\begin{eqnarray*}
\sigma_i\sigma_j&=&\sigma_j\sigma_i,
\quad |i-j|>1,\\
\sigma_i\sigma_{i+1}\sigma_i&=&\sigma_{i+1}\sigma_i\sigma_{i+1},
\end{eqnarray*} and the quadratic relations
$\sigma_i^2=z\sigma_i+1$, with $z=s-s^{-1}$, giving the
alternative form $(\sigma_i-s)(\sigma_i+s^{-1})=0$.

A simple adjustment of the skein relations, as in \cite{AistonMorton},
allows for a skein model  $H_n$ whose parameters can be readily
adapted to match any of the different appearances of the algebra,
\cite{Murphy}.

\begin{Def} Write $H_{n,p}$ for the skein $\CS(F)$ of $(n,p)$-tangles, where
$F$ is the rectangle with $n$ outputs and $p$ inputs at the top, and matching
inputs and outputs at the bottom.
\end{Def}
\begin{center}
$F\ =\ $ \Fnp\\ 
\end{center} 
There is again a natural algebra structure on $H_{n,p}$ induced by
placing  tangles one above the other. When $p=0$ we have the Hecke
algebra $H_n=H_{n,0}$. The resulting algebra $H_{n,p}$ has been studied by
Kosuda and Murakami, \cite{Kosuda}, in the context of $sl(N)_q$ endomorphisms
of the module $V^{\otimes n}\otimes \overline{V}^{\otimes p}$, where $V$ is
the fundamental $N$-dimensional module. Hadji  gives
an explicit skein-theoretic basis for it in \cite{Hadji}, although in
his account the string orientations are the reverse of those used
here.

\subsection{The  annulus}

The Homfly
skein of the annulus, ${\cal C}$, as discussed in
\cite{Turkey} and originally in \cite{Turaev},  is the framed Homfly skein
$\CS(S^1 \x I)$. An element $X\in\CC$ will be indicated schematically as
\[\XC\ .\]
The skein ${\CC}$ has a product induced by
placing one annulus outside another,  under which $\cal C$ becomes a
commutative algebra;
\[\XYC\ =\ \XCYC\ .\]

The closure map $H_{n,p}\to {\cal C}$, induced by taking an $(n,p)$-tangle
$T$ to its closure $\hat {T}$ in the annulus, is defined by \[\hat {T}\ =\
\TC\ .\] This is a $\Lambda$-linear map, whose image we  call ${\CC}_{n,p}$.
Every diagram in the annulus represents an element  in some $\CC_{n,p}$.

Turaev \cite{Turaev} showed  that $\CC$ is freely generated as an algebra
by the set
$\{A_m:m\in\ints\}$ where  $A_m, m\ne 0$, is  represented
by the closure of the braid $\si_{|m|-1}\cdots\si_2\si_1$,
\[A_m\ =\ \Am\ .\] 
The orientation of the curve around  the annulus is   counter-clockwise for
positive
$m$ and clockwise for negative $m$.  The element $A_0$ is the identity
element and is represented by the empty diagram.

The algebra $\CC$ is the product of two subalgebras
$\CC^+$ and $\CC^-$  generated by $\{A_m:m\in\ints,m\ge0\}$ and
$\{A_m:m\in\ints,m\le0\}$ respectively.

The algebra $\CC^+$ is spanned by the subspaces $\CC_{n,0}$ which arise from
the closure of elements in $H_n$, and has been studied quite extensively
\cite{AistonMorton, Lukac, Murphy, Kawagoe}. In
particular, there is a good basis of $\CC^+$ consisting of closures of
certain idempotents of $H_n$.

The basis elements $\{\qlm\}$ for the whole of $\CC$, which we construct
here,  are
identified in section \ref{sect3} as  eigenvectors of the meridian
maps 
defined in section \ref{sect1.7}, following the Schur function methods of
Lukac \cite{Lukac} for
$\CC^+$.  This fact suggests that it should
be useful and relatively easy to express other elements of $\CC$ in
this basis.  Lukac also
shows that the Gyoja-Aiston idempotents for $H_n$,
\cite{AistonMorton}, close to the basis elements $Q_{\la,\phi}$ which
span $\CC^+$. Our results here complement the work of Kosuda and
Murakami \cite{Kosuda} on skein-based  idempotents
in
$H_{n,p}$ and their closures in $\CC$ along the same lines.

\subsection{The annulus with two boundary points}

Denote by $\CA$ the Homfly skein  of the annulus with one input and one
output  boundary point, 
\[ \Annulus\ \]
one
on each boundary component, as indicated.

In a similar way to $\CC$, the skein $\CA$ becomes an algebra under the
product induced by placing one annulus outside another.  The identity element
$e\in\CA$ is represented by the diagram \[e\ =\ \Aid\ \]
obtained by joining the two boundary points by a single straight arc.

A further element $a$ of $\CA$ is  represented by \[a\ =\ \Aa\ .\]  Powers
of this element,
$a^m$ for $m\in\ints$, are also represented by diagrams with no crossings.
For example, 
\[a^{-1}\ =\ \Aainv\ ,\ a^2\ =\ \Aaa\ .\]

As an algebra, $\CA$
 is commutative.  A straightforward skein theoretic proof is given in
\cite{Murphy}, although unlike the case of
$\CC$ the result is not immediately obvious.

\begin{remark}
A skein which is linearly isomorphic to $\CA$ is used by Kawagoe
\cite{Kawagoe} and other authors.  Their version is based on the annulus with
input and output points on the same component.  More recently $\CA$ has come
into use, as its  algebraic properties allow for some satisfyingly
clean proofs.  See also \cite{Murphy}, and work by Lukac
\cite{Lukac}.

\end{remark}

\subsection{Skein involutions}

We  make use of two readily defined involutions on the skeins $\CC$
and
$\CA$ in our subsequent calculations.

\begin{itemize}
\item
For every surface $F$ we can define the {\em mirror map},
$\mirror:\CS(F)\rightarrow\CS(F)$ as follows.
  For a tangle $T$ in $F$
define $\bar{T}$ to be $T$ with all its crossings switched. In the
coefficient ring  $\La$ define conjugation by $\bar{s}=s^{-1}$ and
$\bar{v}=v^{-1}$. The operation of switching crossings in tangles and
conjugating coefficients respects the skein relations and so induces a
conjugate-linear automorphism $\mirror:\CS(F)\rightarrow\CS(F)$, which we
call the { mirror map} on $\CS(F)$.

\item
 Rotation of
diagrams in the annulus $S^1\x I$ by $\pi$ about the horizontal axis
through the distinguished boundary points induces a linear automorphism of
each of  the two skeins $\CC$ and $\CA$ of the annulus which we denote by
$*$ in each case. 
\end{itemize}

Thus for the element $a\in \CA$ we have $a^*=a^{-1}$, while
$\bar{a}=a$. We can  see that $(A_m)^*=A_{-m}$, so that
$(\CC^+)^*=\CC^-$. It is also readily verified that
$(\CC_{n,p})^*=\CC_{p,n}$, and that $\CC_{n,p}$ is invariant under the
mirror map, since $H_{n,p}$ is.

\subsection{Some basic elements in the annulus}

 In $H_n$ there is a known set of idempotent elements, $E_\la$,
one for each partition $\la$ of $n$.  These elements were originally given a
purely algebraic description by Gyoja \cite{Gyoja}, and were subsequently
given a skein picture based on the Young diagram for $\la$ by the second
author and Aiston \cite{AistonMorton}. The closure of $E_\la$ gives the
element $Q_\la \in \CC^+$, making up the basis elements used in
\cite{AistonMorton}, and indirectly those used in \cite{Lukac}. Our
construction makes use only of the sequence of elements $\{h_n\}\in\CC^+$
which arise from the single row Young diagrams, and the corresponding
elements $\{h_n^*\}\in\CC^-$.

The two simplest
idempotents in $H_n$ correspond to the single row and single column Young
diagrams, as described in \cite{Turkey}. In the interpretation of the
subalgebra $\CC^+$ of $\CC$ as the algebra of symmetric functions in
infinitely many variables the closures of these idempotents correspond to
the complete symmetric function and the elementary symmetric function of
degree $n$ respectively. We give  brief details of the definition and
properties of these idempotents and their closure.

Let $w_\pi$ be the positive permutation braid (\cite{EM94})
corresponding to
$\pi\in S_n$.  Define two quasi-idempotents $a_n, b_n \in H_n$ by
\[
a_n=\sum_{\pi\in
S_n}s^{l(\pi)}w_\pi\qquad\mbox{and}\qquad b_n=\sum_{\pi\in
S_n}(-s)^{-l(\pi)}w_\pi,
\]
where $l(\pi)=\text{wr}(w_\pi)$, the writhe of the braid $w_\pi$.

Write $g_n=1+s\si_{n-1}+s^2\si_{n-1}\si_{n-2}+\ldots
+s^{n-1}\si_{n-1}\si_{n-2}\ldots\si_1$, where $\{\si_i\}$ are the usual
elementary braid generators of the braid group $B_n$.

We have $g_{n+1}=1+s\si_n g_n$, and also the immediate skein relation
\beqnn\Gnplus \quad = \quad \Gshift\quad  + s^n \quad \Gcross
\label{eqn:gn}
\eeqnn
for tangles on $n+1$ strings.

There is an algebra homomorphism $\ph_s:H_n\to \La$ defined on the
elementary braids $\sigma_i$ by $\ph_s(\sigma_i)=s$.

The following two lemmas are readily established, as in \cite{Turkey}: 
\begin{lemma}
 $a_n=a_{n-1}g_n$.
\end{lemma}

\begin{lemma}
For any $g\in H_n$ we have $a_ng=\ph_s(g)a_n=g a_n$.
\end{lemma}

The element $a_n$ then satisfies
\[  a_n^2=\ph_s(a_n)a_n=\ph_s(a_{n-1})\ph_s(g_n)a_n.  \]

Now since $\ph_s(g_n)=1+s^2+\ldots +s^{2n-2}=s^{n-1}[n]$ with
$[k]=\frac{s^k-s^{-k}}{s-s^{-1}}$, we have the immediate
corollary:
\begin{corollary} \label{hn}
We can write 
\[ s^{n-1}[n]h_n=h_{n-1}g_n, \]
where $h_n=a_n/\ph_s(a_n)$ is the true idempotent.
\end{corollary}

The element $h_n$ constructed above is the idempotent which corresponds to
the single row Young diagram with $n$ cells.  The single column idempotent,
denoted $e_n$, is constructed in an analogous way from $b_n$.  It can be
obtained from $h_n$ by using $-s^{-1}$ in place of $s$.

With a slight abuse of the notation we shall write $h_n,e_n\in\CC$ for the
closures $\hat{h}_n$, $\hat{e}_n$ of these idempotents as elements of
$\CC$.

The skein $\CC^+$ is spanned by the monomials in $\{h_m:m\ge 0\}$.
  The skein
$\CC^-$ is thus spanned by monomials in $\{h_l^*:l\ge 0\}$.
The whole skein $\CC$ is then spanned by monomials in
$\{h_l^*, h_m:l,m\ge 0\}$.

\subsection{The meridian maps of the skein $\CC$} \label{sect1.7}
We introduce here the linear endomorphisms $\ph,\bar{\ph}$ of $\CC$ which
are of central importance in this work.

Take a diagram $X$ in the annulus and link it once with a simple meridian
loop, oriented in either direction, to give diagrams $\ph(X)$ and
$\bar{\ph}(X)$ in the annulus as indicated. This induces linear endomorphisms
$\ph,\bar{\ph}$ of the skein $\CC$, called the {\em meridian maps}.

\[
\begin{array}{ccc}
\varphi:\CC&\rightarrow&\CC \\
\Xone&\mapsto&
\Xtwo
\end{array} ,\ \text{ and }\quad
\begin{array}{ccc}
\bar{\varphi}:\CC&\rightarrow&\CC \\
\Xone&\mapsto&
\Xthree
\end{array}.
\]

Each subspace $\CC_{n,p}$ is invariant under $\ph$ and $\bar\ph$. In
\cite{Lukac} the eigenvectors of $\ph$ on $\CC_{n,0}$ are identified
with the closures $Q_\la$ of the idempotents in $H_n$, for partitions
$\la$ of
$n$, and are also expressed as explicit integer polynomials in
$\{h_r\}$, as in
\cite{Kawagoe}. Lukac makes use of the skein $\CA$ in doing this, and our
methods here are an extension of his work.

In \cite{MH02} we calculated the complete set of eigenvalues
$\{t_{\la,\mu}\}$ of $\ph$, where $(\la,\mu)$ run over all pairs of
partitions of integers, and showed that they are all distinct.  Our
goal in this paper  is an
explicit expression, as an integer polynomial  in $\{h_n\}$ and
$\{h_n^*\}$, for an eigenvector
$\qlm$ corresponding to each eigenvalue $t_{\mu,\la}$. The
general result is given in section \ref{sect3.2}, along with an
explicit example where $\la$ and $\mu$ are the partitions with
parts $4,2,2$ and $3,2$ respectively.

Under the mirror map and the involution $*$ it is clear that
\beqnn
\overline{\ph(X)}&=&\bar{\ph}(\bar{X})\\
(\ph(X))^*&=&\bar{\ph}(X^*).
\eeqnn

It can be established  that $h_n$ is fixed by the mirror map, and is
an eigenvector of $\ph$, see for example \cite{Lukac}. Use of the mirror map
and
$*$ then show at once that $h_n$ and $h_n^*$ are eigenvectors both of $\ph$
and of
$\bar{\ph}$. The eigenvalues are listed explicitly in the following lemma.
\begin{lemma}\label{l1.4} 
\beqn
\phia(h_n) &=& (v^{-1}(s^{2n-1}-s^{-1})+\de)h_n, \\
\phia(h_n^*) &=& (v(s^{-2n+1}-s)+\de)h_n^*, \\
\phib(h_n) &=& (v(s^{-2n+1}-s)+\de)h_n, \\
\phib(h_n^*) &=& (v^{-1}(s^{2n-1}-s^{-1})+\de)h_n^*.
\eeqn
\end{lemma}

\subsection{Relations between $\CC$ and $\CA$} \label{cca}

Two algebra homomorphisms $l,r:\CC\to \CA$ can be induced by  placing an
element of $\CC$ respectively under or over the simple arc $e$ joining the
distinguished boundary points of the annulus. Thus
\[
l(X)=\LX \qquad
\mbox{and}\qquad r(X)=
\RX\ .
\] 
The effects of the mirror
map and the involution $*$ are readily seen to give
\beqnn
\overline{l(X)}&=&r(\overline{X})\\
(l(X))^*&=&r(X^*).
\eeqnn
Then $\overline{l(h_n)}=r(h_n)$, where $h_n\in\CC$ is the closure of the
idempotent in $H_n$.

 Define a  closure map $\diamond:\CA \rightarrow\CC$ on diagrams  by
joining the distinguished boundary points with a simple arc lying   above
the annulus.   This map interacts with $l$ and $r$ to give
\beqn
\diamond(r(X))&=&\bar{\ph}(X),\\
\diamond(l(X))&=&\de X,\\
\diamond(Yl(X))&=&\diamond(Y)X,
\eeqn
for any $X\in\CC, Y\in \CA$.

We shall construct eigenvectors of $\ph$ and $\bar{\ph}$ as determinants
of matrices with entries $h_n, h_n^*$ in $\CC$. When $M$ is any such a
matrix, with
$X=\det M$, we have $\bar{\ph}(X)=\diamond(r(\det M))=\diamond(\det
r(M))$, where $r(M)$ is the matrix with entries in $\CA$ given by applying
the algebra homomorphism $r$ to all entries in $M$. Our main calculation
in the next section involves  determinants of  matrices such as $r(M)$ with
entries in
$\CA$. In advance of this
we introduce some relations in $\CA$ which  allow us to perform certain
column operations on the matrices.

Define elements $y_n\in\CA$, fixed under the mirror map,  by
\[
y_n:=[n+1]\x\Yn
\]
where $h_{n+1}$ here is the idempotent in $H_{n+1}$.
\begin{lemma} We have the relation
\beqnn\label{eqn:yn1}
y_n=s^{-1}ay_{n-1}+l(h_n).
\eeqnn
\end{lemma}

\begin{proof}
Apply Corollary \ref{hn} to write $[n+1]h_{n+1}=s^{-n}h_n g_{n+1}$. Use the
skein relation (\ref{eqn:gn}) to write $g_{n+1}$, and hence $y_n$, as the
sum of two parts. One of these gives $l(h_n)$ at once and the other gives
$s^{-1}ay_{n-1}$ on moving the box $g_n$ in (\ref{eqn:gn}) round the annulus
to combine with $h_n$ from above.
\end{proof}

Apply the mirror map to
 (\ref{eqn:yn1}), noting that  $\bar{s}=s^{-1}$, to
get 
\beqnn\label{eqn:yn2}
y_n=say_{n-1}+r(h_n).
\eeqnn

Since $a^*=a^{-1}$ we can apply $*$ to equations (\ref{eqn:yn1}) and
(\ref{eqn:yn2}) to obtain

\beqnn\label{eqn:yn*1}
y_n^*&=&s^{-1}a^{-1}y_{n-1}^*+r(h_n^*)\\
\label{eqn:yn*2}
y_n^*&=&sa^{-1}y_{n-1}^*+l(h_n^*).
\eeqnn

We rewrite  equations (\ref{eqn:yn1}) and (\ref{eqn:yn*2}) as
\beqnn
y_n&=&s^{-1}ay_{n-1}+l_n, \label{eqn:yn1a}\\
y_{n-1}^*&=&s^{-1}ay_n^*+L_{n-1}, \label{eqn:yn*3} 
\eeqnn
with $l_n=l(h_n)$
and $L_{n-1}=-s^{-1}al(h_n^*)$, for ease of use in our later matrix work.

Similarly we rewrite (\ref{eqn:yn2}) and (\ref{eqn:yn*1}) as
\beqnn 
y_n&=&say_{n-1}+r_n \label{eqn:yn2a}\\
y_{n-1}^*&=&say_n^*+R_{n-1} \label{eqn:yn*4}.
\eeqnn
with $r_n=r(h_n)$ and $R_{n-1}=-sar(h_n^*)$.

From equations (\ref{eqn:yn1a}) and (\ref{eqn:yn2a})  we get 
\beqnn
(s-s^{-1})y_n&=&sl_n-s^{-1}r_n ,\label{eqn:yna} 
\eeqnn
and from
(\ref{eqn:yn*3}) and (\ref{eqn:yn*4}) we get
\beqnn
(s-s^{-1})y_{n-1}^*&=&sL_{n-1}-s^{-1}R_{n-1}.
\label{eqn:ynb}
\eeqnn

\section{Matrix-based calculations}

In this section we  work with square matrices whose entries, in
$\CC$ or
$\CA$, are arranged in a fairly restricted format. Because of the
restrictions we can use a simplified notation for them, which we
now describe.

\subsection{Notation}

The entries in each matrix are drawn from sequences such as
$h=\{h_n\}, y^*=\{y_n^*\}, R=\{R_n\}$ with integer subscripts.  On any
given row the subscripts  either increase in steps of 1 along the
row, (a {\em standard} row), or decrease in steps of 1, (a {\em
starred} row). The subscript for every entry of the matrix is then
determined by knowing which rows are starred, and the subscripts 
to be used in the first column.

 We shall always use matrices in which the starred rows are the first
$k^*$ rows, followed by $k$ standard rows. The column vector $\bf v$ of
subscripts to be used for the first column will be called the {\em index
vector} for the matrix.

The matrices used  have the further restriction that in any given column
the entries in the starred rows are drawn from some fixed sequence,
$h,R$, for example, while all entries in the standard rows of a given
column  again come from some fixed sequence.

By way of example the $8\x8$ matrix $M$ below has entries in such a format.

\[M=\begin{pmatrix}
a_3^* & a_2^* & a_1^* & b_0 & b_{-1} & c_{-2}^* & c_{-3}^* & c_{-4}^*\\ 
a_5^*&a_4^*&a_3^*&b_2&b_1&c_0^* & c_{-1}^* & c_{-2}^* \\
a_4^* & a_3^* & a_2^* & b_1 & b_0 & c_{-1}^* & c_{-2}^* & c_{-3}^*\\
a_{-1}&a_0&a_1&d_2&c_3&c_4&c_5&c_6 \\ 
a_1 & a_2 & a_3 & d_4 & c_5 & c_6 & c_7 & c_8 \\
a_2 & a_3 & a_4 & d_5 & c_6 & c_7 & c_8 & c_9 \\
a_1 & a_2 & a_3 & d_4 & c_5 & c_6 & c_7 & c_8 \\
a_0&a_1&a_2&d_3&c_4 & c_5 & c_6 & c_7\\
\end{pmatrix}.\]
In $M$ the first 3 rows are starred and the index vector is \[{\bf
v}=\begin{pmatrix}3\\
5\\
4\\
\hline
-1\\
1\\
2\\
1\\
0
\end{pmatrix}
\]

The complete matrix $M$ can be recovered, given $\bf v$, once the sequences
to be used in the starred and standard rows are known.

This information is given by the {\em template matrix}, in this case
\[\begin{pmatrix}a^*&a^*&a^*&b&b&c^*&c^*&c^*\\
a&a&a&d&c&c&c&c
\end{pmatrix}.\] In general the template matrix is the $2\x(k^*+k)$ matrix
giving the sequences from which the starred and standard
entries in each column  are drawn.

As a further notational simplification we use the symbol $\dashv$ in the
template matrix to indicate that the same sequence is repeated, and the
symbol 
$\cdots$ where the choice of sequence is unimportant. We label a column or
columns in the template as necessary to show where a repetition ends.

With this notation the template matrix for $M$ can  be written

\[\hspace{-8pt}\bordermatrix{&&&\overset{4}{\scriptscriptstyle{\downarrow}}&&&\cr
&a^*&\dashv&b&b&c^*&\dashv\cr
&a&\dashv&d&c&&\dashv}\]
and the complete matrix $M$ with index vector $\bf v$ can then denoted simply
by 
\[M=\hspace{-8pt}\bordermatrix{&&&\overset{4}{\scriptscriptstyle{\downarrow}}&&&\cr
&a^*&\dashv&b&b&c^*&\dashv\cr
&a&\dashv&d&c&&\dashv}_{\bf v}.\]

For most of our calculations the index vector $\bf v$ will be fixed, but the
choice of $\bf v$ will be unimportant, and we  frequently suppress it,
giving simply the template.
\subsection{Calculations}

Throughout this section we  use matrices whose entries are drawn
from the sequences of elements of $\CA$ or $\CC$ described in 
section \ref{cca}. The sequences will be denoted by the  letters used
there, $h,h^*,y,y^*,l,r,L,R$. The definitions can in all cases be
extended to allow for negative index by setting
$h_j=h^*_j=y_j= y^*_j= l_j= r_j= L_j= R_j=0$ when $j<0$.
The equations at the end of section \ref{cca} relating the sequences
continue to hold in the extended range.
\begin{lemma}
For each $i\ge1$ and fixed index vector we have the following equations
in $\CA$.
\beqnn
\det\hspace{-8pt}\bordermatrix{&&\overset{i}{\scriptscriptstyle{\downarrow}}&&\cr
&\cdots&y^*&R&\cdots\cr
&\cdots&y&r&\cdots}&=&\det\hspace{-8pt}\bordermatrix{&&\overset{i}{\scriptscriptstyle{\downarrow}}&&\cr
&\cdots&y^*&y^*&\cdots\cr
&\cdots&y&y&\cdots} \label{eqn:yr}\\
\det\hspace{-8pt}\bordermatrix{&&\overset{i}{\scriptscriptstyle{\downarrow}}&&\cr
&\cdots&y^*&L&\cdots\cr
&\cdots&y&l&\cdots}&=&\det\hspace{-8pt}\bordermatrix{&&\overset{i}{\scriptscriptstyle{\downarrow}}&&\cr
&\cdots&y^*&y^*&\cdots\cr
&\cdots&y&y&\cdots} \label{eqn:yl}
\eeqnn
\end{lemma}
\begin{proof}
The column operation $ C_{i+1}\mapsto saC_i+C_{i+1} $ applied to the
matrix on the lefthand side of equation (\ref{eqn:yr}) gives the matrix on
the righthand side,
 using (\ref{eqn:yn*4}) on starred rows and (\ref{eqn:yn2a}) on standard
 rows.

 The column operation $
C_{i+1}\mapsto s^{-1}aC_i+C_{i+1} $ works similarly for equation
(\ref{eqn:yl})
 using (\ref{eqn:yn*3}) on starred rows and (\ref{eqn:yn1a}) on standard
 rows.
\end{proof}

\begin{corollary} \label{yrl}
For each $j\ge 1$ we have \beqn
\det\hspace{-8pt}\bordermatrix{&&\overset{j}{\scriptscriptstyle{\downarrow}}&&\cr
&\cdots&y^*&R&\dashv\cr
&\cdots&y&r&\dashv}=
\det\hspace{-8pt}\bordermatrix{&&\overset{j}{\scriptscriptstyle{\downarrow}}&\cr
&\cdots&y^*&\dashv\cr
&\cdots&y&\dashv}
=\det\hspace{-8pt}\bordermatrix{&&\overset{j}{\scriptscriptstyle{\downarrow}}&&\cr
&\cdots&y^*&L&\dashv\cr &\cdots&y&l&\dashv}
\eeqn
\end{corollary}
\begin{proof}
Apply equation (\ref{eqn:yr}) repeatedly, for $i$ from $j$ onwards, to get
the first equation, and similarly use (\ref{eqn:yl}) repeatedly to get the
second equation.
\end{proof}

\begin{lemma} \label{y} For each $j\ge1$ we have
\beqn(s-s^{-1})\det\hspace{-8pt}\bordermatrix{&&\overset{j}{\scriptscriptstyle{\downarrow}}&\cr
&\cdots&y^*&\cdots\cr
&\cdots&y&\cdots}&=&
s\det\hspace{-8pt}\bordermatrix{&&\overset{j}{\scriptscriptstyle{\downarrow}}&\cr
&\cdots&L&\cdots\cr
&\cdots&l&\cdots}\\
&&-s^{-1}\det\hspace{-8pt}\bordermatrix{&&\overset{j}{\scriptscriptstyle{\downarrow}}&\cr
&\cdots&R&\cdots\cr
&\cdots&r&\cdots}.
\eeqn
\end{lemma}
\begin{proof} Expand the determinants by the $j$th column, using equations
(\ref{eqn:ynb}) on the starred rows and (\ref{eqn:yna}) on the standard rows.
\end{proof}

\begin{notation}
Write
\beqn
\De_j&=&
\det\hspace{-8pt}\bordermatrix{&&&\overset{j+1}{\scriptscriptstyle{\downarrow}}&\cr
&L&\dashv&R&\dashv\cr
&l&\dashv&r&\dashv},\ 0\le j<k^*+k,\\
\De_{k^*+k}&=&\det\hspace{-8pt}\bordermatrix{&&\cr
&L&\dashv\cr &l&\dashv},
\eeqn
for any choice of index vector.
\end{notation}

\begin{lemma}\label{delta}
We have
\beqn
(s-s^{-1})\det\hspace{-8pt}\bordermatrix{&&&\overset{j}{\scriptscriptstyle{\downarrow}}&&\cr
&L&\dashv&y^*&L&\dashv\cr
&l&\dashv&y&l&\dashv} =s\De_j-s^{-1}\De_{j-1},\ 1\le j\le k^*+k.
\eeqn
\end{lemma}
\begin{proof}
\beqn
(s-s^{-1})\det\hspace{-8pt}\bordermatrix{&&&\overset{j}{\scriptscriptstyle{\downarrow}}&&\cr
&L&\dashv&y^*&L&\dashv\cr
&l&\dashv&y&l&\dashv}&=&
(s-s^{-1})\det\hspace{-8pt}\bordermatrix{&&&\overset{j}{\scriptscriptstyle{\downarrow}}&&\cr
&L&\dashv&y^*&R&\dashv\cr
&l&\dashv&y&r&\dashv},\\
&& \text{ by corollary \ref{yrl}}\\
&=&s\De_j-s^{-1}\De_{j-1}, \text{ by lemma \ref{y}}
\eeqn
\end{proof}

Expanding further by lemma \ref{y} gives
\beqnn
s\De_j-s^{-1}\De_{j-1}=s \det\hspace{-8pt}\bordermatrix{&&\cr
&L&\dashv\cr
&l&\dashv}
-s^{-1}\det\hspace{-8pt}\bordermatrix{&&&\overset{j}{\scriptscriptstyle{\downarrow}}&&\cr
&L&\dashv&R&L&\dashv\cr &l&\dashv&r&l&\dashv}.
\label{eqn:dj}
\eeqnn

\begin{corollary}\label{LRL} We have
\beqn
 \sum_{j=1}^{k^*+k-1}s^{2j}
\det\hspace{-8pt}\bordermatrix{&&\cr
&L&\dashv\cr
&l&\dashv}
+\det\hspace{-8pt}\bordermatrix{&&\cr
&R&\dashv\cr
&r&\dashv}\\
=\sum_{j=1}^{k+k^*}s^{2j-2}
\det\hspace{-8pt}\bordermatrix{&&&\overset{j}{\scriptscriptstyle{\downarrow}}&&\cr
&L&\dashv&R&L&\dashv\cr &l&\dashv&r&l&\dashv}.
\eeqn
\end{corollary}
\begin{proof}
Multiply both sides of (\ref{eqn:dj}) by $s^{2j-1}$ and sum to get
\beqn
s^{2(k^*+k)}\De_{k^*+k}-\De_0&=&\sum_{j=1}^{k+k^*}s^{2j}
\det\hspace{-8pt}\bordermatrix{&&\cr
&L&\dashv\cr
&l&\dashv}\\
&&-
\sum_{j=1}^{k+k^*}s^{2j-2}
\det\hspace{-8pt}\bordermatrix{&&&\overset{j}{\scriptscriptstyle{\downarrow}}&&\cr
&L&\dashv&R&L&\dashv\cr &l&\dashv&r&l&\dashv},
\eeqn
and then rearrange the terms.
\end{proof}

\begin{lemma}\label{detlh} For any index vector $\bf w$ we have
\beqn
 s^{-k^*}\left(\sum_{j=1}^{k^*+k-1}s^{2j}\right)
\det\hspace{-8pt}\bordermatrix{&&\cr
&l(h^*)&\dashv\cr
&l(h)&\dashv}_{\bf w}
+s^{k^*}\det\hspace{-8pt}\bordermatrix{&&\cr
&r(h^*)&\dashv\cr
&r(h)&\dashv}_{\bf w}\\
=s^{-k^*}\sum_{j=1}^{k+k^*}
\det\hspace{-8pt}\bordermatrix{&&&\overset{j}{\scriptscriptstyle{\downarrow}}&&\cr
&l(h^*)&\dashv&s^{2j} r(h^*)&l(h^*)&\dashv\cr
&l(h)&\dashv&s^{2j-2}r(h)&l(h)&\dashv}_{\bf w}.
\eeqn
\end{lemma}

\begin{proof}
Recall that $L_{n-1}=(-a)s^{-1}l(h_n^*)$ and $R_{n-1}=(-a)sr(h_n^*)$. 
Use corollary \ref{LRL} with index vector $\bf v$  given by reducing all
starred indices of $\bf w$  by 1, and leaving the standard indices unaltered.
Then  take
out a factor of $-a$ from each starred row to get
\beqn
 (-a)^{k^*}\left(\sum_{j=1}^{k^*+k-1}s^{2j}\right)
\det\hspace{-8pt}\bordermatrix{&&\cr
&s^{-1}l(h^*)&\dashv\cr
&l(h)&\dashv}_{\bf w}
+(-a)^{k^*}\det\hspace{-8pt}\bordermatrix{&&\cr
&sr(h^*)&\dashv\cr
&r(h)&\dashv}_{\bf w}\\
=(-a)^{k^*}\sum_{j=1}^{k+k^*}
s^{2j-2}\det\hspace{-8pt}\bordermatrix{&&&\overset{j}{\scriptscriptstyle{\downarrow}}&&\cr
&s^{-1}l(h^*)&\dashv&s r(h^*)&s^{-1}l(h^*)&\dashv\cr
&l(h)&\dashv&r(h)&l(h)&\dashv}_{\bf w}.
\eeqn

The factor $(-a)^{k^*}$ can be cancelled. Extract $k^*$ further factors of
$s$ or $s^{-1}$ from the starred rows, and insert the factor $s^{2j-2}$ into
the $j$th column of the appropriate matrix to complete the proof.
\end{proof}

\section{Eigenvectors of the meridian maps} \label{sect3}
We now come to the main result, giving potential eigenvectors $A_{\bf w}$ in
$\CC$ for the meridian maps. We show later that sufficiently many of these
vectors are non-zero, to provide our claimed basis $\qlm$ by suitable
choices of index vector $\bf w$.
\begin{theorem} \label{t2.7}
Let $\bf w$ be any index vector. Then \[A_{\bf w}=\det\hspace{-8pt}\bordermatrix{&&\cr
&h^*&\dashv\cr
&h&\dashv}_{\bf w} \in \CC\]
satisfies $\bar{\ph}(A_{\bf w})=c_{\bf w}A_{\bf w}$, and hence also
${\ph}(A_{\bf w})=\overline{c_{\bf w}}A_{\bf w}$, for some
$c_{\bf w}\in\La$.
\end{theorem}

\begin{proof} Write $M_{\bf w}=\hspace{-8pt}\bordermatrix{&&\cr
&h^*&\dashv\cr
&h&\dashv}_{\bf w}$, with entries $m_{i\,j}$ and cofactors
$M_{i\,j}\in\CC$.
Write $c_{i\,j}$ for the entries in the distinguished $j$th column of
\[
\hspace{-8pt}\bordermatrix{&&&\overset{j}{\scriptscriptstyle{\downarrow}}&&\cr
&l(h^*)&\dashv&s^{2j} r(h^*)&l(h^*)&\dashv\cr
&l(h)&\dashv&s^{2j-2}r(h)&l(h)&\dashv}_{\bf w}.
\]

The determinant of this matrix, expanded by the $j$th column is
\[\sum_{i=1}^{k^*+k}c_{i\,j}l(M_{i\,j}).\]
Consequently  lemma \ref{detlh}
gives
\beqnn
s^{-2k^*}\left(\sum_{j=1}^{k^*+k-1}s^{2j}\right)\det l(M_{\bf w})+\det
r(M_{\bf w})=s^{-2k^*}\sum_{i,j=1}^{k^*+k}c_{i\,j}l(M_{i\,j}).
\label{eqn:lraw}
\eeqnn

Now apply the closure map $\diamond$ to the elements of $\CA$ on each
side of equation (\ref{eqn:lraw}). Recall that 
\beqn
\diamond(l(X))&=&\de X,\\
\diamond(r(X))&=&\bar{\ph}(X),\\
\diamond(Yl(X))&=&\diamond(Y) X, \text{ for any }X\in\CC,Y\in\CA,
\eeqn
and that $r$ and $l$ are algebra homomorphisms, so that $\det
l(M_{\bf w})=l(A_{\bf w})$ and $\det r(M_{\bf w})=r(A_{\bf w})$.
We then have
\beqnn
s^{-2k^*}\left(\sum_{j=1}^{k^*+k-1}s^{2j}\right)\de A_{\bf
w}+\bar{\ph}A_{\bf
w}=s^{-2k^*}\sum_{i,j=1}^{k^*+k}\diamond(c_{i\,j})M_{i\,j}.
\label{eqn:aw}
\eeqnn

We now show that $\diamond(c_{i\,j})$ is a multiple of $m_{i\,j}$, and that the
righthand side of equation (\ref{eqn:aw}) can be written as a multiple of
$A_{\bf w}$, thus proving the theorem.

We deduce this from the following  lemma.
\begin{lemma} \label{cij} There exist $\alpha_i,\be_j\in\CA,1\le i,j\le
k^*+k$, such that
$$\diamond(c_{i\,j})=(\alpha_i+\be_j)m_{i\,j}.$$
\end{lemma}
\begin{corollary} 
\beqn
\sum_{i,j=1}^{k^*+k}\diamond(c_{i\,j})M_{i\,j}
&=&\sum_{i,j=1}^{k^*+k}(\alpha_i+\be_j)m_{i\,j}M_{i\,j}\\
&=&\sum_{i,j=1}^{k^*+k}\alpha_i
m_{i\,j}M_{i\,j}+\sum_{i,j=1}^{k^*+k}\be_j m_{i\,j}M_{i\,j}\\
&=&\sum_{i=1}^{k^*+k}\alpha_i\left(\sum_{j=1}^{k^*+k}m_{i\,j}M_{i\,j}\right)+
\sum_{j=1}^{k^*+k}\be_j\left(\sum_{i=1}^{k^*+k}m_{i\,j}M_{i\,j}\right)\\
&=&\sum_{i=1}^{k^*+k}\alpha_i A_{\bf w}
+\sum_{j=1}^{k^*+k}\be_j A_{\bf w},
\eeqn
using the expansion of $A_{\bf w}$ by the $i$th row in the first sum and by
the
$j$th column in the second.
\end{corollary}

\begin{proof}[Proof of lemma \ref{cij}]
Since \beqn c_{i\,j}=\begin{cases}
s^{2j}r(m_{i\,j}) & \text{if $i\le k^*$,}\\
s^{2j-2}r(m_{i\,j}) & \text{if $i> k^*$,}
\end{cases}
\eeqn we have 
\beqn \diamond( c_{i\,j})=\begin{cases}
s^{2j}\bar{\ph}(m_{i\,j}) & \text{if $i\le k^*$,}\\
s^{2j-2}\bar{\ph}(m_{i\,j}) & \text{if $i> k^*$.}
\end{cases}
\eeqn

In every case $m_{i\,j}$ is either $h_n^*$ or $h_n$ for some $n$ depending
on the index vector $\bf w$, and is hence an eigenvector of $\bar{\ph}$. 

For $i\le k^*$ we have $m_{i\,j}=h_n^*$ with $n=w_i -j+1$ and
$\bar{\ph}(m_{i\,j})=(v^{-1}s^{2n-1}+\de-v^{-1}s^{-1})m_{i\,j}$, while for
$i>k^*$ we have $m_{i\,j}=h_n$ with $n=w_i +j-1$ and
$\bar{\ph}(m_{i\,j})=(v^{-1}s^{-2n+1}+\de-vs)m_{i\,j}$.
 In each case
$\diamond(c_{i\,j})=(\alpha_i+\be_j)
m_{i\,j}$, where
 \beqn
\alpha_i&=&\begin{cases} v^{-1}s^{2w_i+1} & \text{if $i\le k^*$,}\\
vs^{1-2w_i} & \text{if $i> k^*$}
\end{cases}\\
\be_j&=&s^{2j}(\de-v^{-1}s^{-1})=s^{2j-2}(\de-vs)\\
&=&\frac{v^{-1}s^{2j-2}-vs^{2j}}{s-s^{-1}}.
\eeqn
\end{proof}  

The theorem is now established, with 
\beqnn c_{\bf
w}&=&s^{-2k^*}\left(\sum_{i=1}^{k^*+k}\alpha_i+\sum_{j=1}^{k^*+k}\be_j
-\de\sum_{j=1}^{k^*+k-1}s^{2j}\right)
\label{eqn:cw}\eeqnn 
\end{proof}

\subsection{Construction of $\qlm$} \label{sect3.2}
Given two partitions $\la$ and $\mu$ with $k$ and $k^*$ parts
respectively we define $\qlm=A_{\bf w}$ for a suitably chosen index
vector $\bf w$ with $k$ standard and $k^*$ starred rows. We choose
$\bf w$ so that the subscripts of the {\em diagonal} entries in the
standard rows are the parts $\la_1,\la_2,\ldots,\la_k$ of $\la$ in
order, while the subscripts  of the {\em diagonal} entries in the
starred rows are the parts $\mu_1,\mu_2,\ldots,\mu_{k^*}$ of $\mu$ in
{\em reverse} order.

The entries in the index vector are then explicitly \[
w_i=\begin{cases}
\mu_{k^*-i+1}+i-1&i\le k^*,\\
\la_{i-k^*}-i+1&i>k^*.
\end{cases}\]

\begin{remark}
Similar determinants are used by Koike \cite{Koike} in giving universal
formulae for the irreducible characters of rational representations of
$GL(N)$, along with interpretations in terms of skew Schur functions.
\end{remark}

\begin{example} 
For the partitions $\la$ and $\mu$ with parts $4,2,2$ and $3,2$
respectively, represented by the Young diagrams
$\la=\fourtwotwo$ and $\mu=\threetwo$ this gives
\[\qlm=\det\begin{pmatrix} {\bf h^*_2}&h^*_1&1&0&0\\
h^*_4&{\bf h^*_3}&h^*_2&h^*_1&1\\
h_2&h_3&{\bf h_4}&h_5&h_6\\
0&1&h_1&{\bf h_2}&h_3\\
0&0&1&h_1&{\bf h_2}
\end{pmatrix}\]
\end{example}
\begin{remark} Unless the matrix $M_{\bf w}$ has two equal rows or an
entirely zero row the
 elements $A_{\bf w}=\det M_{\bf w}$ of $\CC$ are $\pm
\qlm$ for some partitions $\la,\mu$.

 If the index vector
$\bf w$ has a repeated entry in the starred rows, or a repeated entry
in the standard rows then $M_{\bf w}$ has a repeated row. When $w_i<0$ for
any
$i\le k^*$ or $w_i<k^*+k$ for any
$i\ge k^*$  the matrix $M_{\bf w}$ will  have a zero row. Otherwise, by
permuting the rows we can assume that the starred entries in $\bf w$
increase with $i$ and the standard entries decrease with $i$, and that
$w_1\ge0$ and
$w_{k^*+k}\ge -k^*-k$. If $w_1=0$ or $w_{k^*+k}=
-k^*-k$ then we can find an index vector with one fewer starred or
standard entries respectively which determines the same element
$A_{\bf w}$, on expanding the determinant by either the first or the last
row. We may then assume that the starred entries increase strictly, with
$w_1>0$ and that the standard entries decrease strictly, with $w_{k^*+k}>
-k^*-k$. The indices on the diagonal then  determine the partitions $\la$
and $\mu$ for which $A_{\bf w}=\qlm$.
\end{remark}

\begin{definition}
Given two partitions $\la$ and $\mu$, write 
\[s_{\la,\mu}=(s-s^{-1})\left(v^{-1}\sum_{x\in\la}s^{2c(x)}
-v\sum_{x\in\mu}s^{-2c(x)}\right)+\delta,\]
where the sum is taken over cells $x$ in the Young diagram of the
partition, and
$c(x)=j-i$ is the {\em content} of the cell $x$ in row
$i$ and column
$j$ of the Young diagram.
\end{definition}
\begin{remark} In \cite{MH02} the notation $t_{\la,\mu}$ is used, with
$t_{\la,\mu}=s_{\mu,\la}$. It is shown there that the set
$\{t_{\la,\mu}\}$ forms a complete set of eigenvalues of $\ph$, each
occurring with multiplicity $1$.
\end{remark}
\begin{theorem}
The element $\qlm\in\CC$ defined above is an eigenvector of the
meridian map
$\ph$, with eigenvalue $s_{\la,\mu}$.
\end{theorem}
\begin{proof} We know already from theorem \ref{t2.7} that either
$\qlm=0$ or it is an eigenvector of $\ph$. We first identify its
eigenvalue as $s_{\la,\mu}$ by showing 
that
$\overline{\ph}(\qlm)=s_{\mu,\la}\qlm$. This is sufficient, since then
$\ph(\qlm)=\overline{s_{\mu,\la}}\qlm=s_{\la,\mu}\qlm$.
 Finally we establish
that
$\qlm\ne 0$ in $\CC$. 

From the explicit formula  (\ref{eqn:cw}) in theorem \ref{t2.7}
for
$c_{\bf w}$ in terms of the index vector $\bf w$ which defines $\qlm$ we
have
\beqn
c_{\bf w}&=&s^{-2k^*}\left(\sum_{i=1}^{k^*+k}\alpha_i+\sum_{j=1}^{k^*+k}\be_j
-\de\sum_{j=1}^{k^*+k-1}s^{2j}\right)\\
&=&v^{-1}\sum_{i=1}^{k^*}s^{2(w_i-k^*)+1}+v\sum_{i=k^*+1}^{k^*+k}s^{-2(w_i+k^*)+1}
\\
&&+s^{-2k^*}\left(\sum_{j=1}^{k^*+k}\be_j
-\de\sum_{j=1}^{k^*+k-1}s^{2j}\right).\\
\eeqn

Now $\be_j=\frac{v^{-1}}{s-s^{-1}}s^{2j-2}-\frac{v}{s-s^{-1}}s^{2j}$ and
$\de=\frac{v^{-1}}{s-s^{-1}}-\frac{v}{s-s^{-1}}$. Then
\beqn
\sum_{j=1}^{k^*+k}\be_j
-\de\sum_{j=1}^{k^*+k-1}s^{2j}&=&
v^{-1}\frac{1}{s-s^{-1}}-v\frac{s^{2k^*+2k}}{s-s^{-1}}.
\eeqn

For our choice of $\bf w$ in terms of $\la$ and $\mu$ we have
\[
\sum_{i=1}^{k^*}s^{2(w_i-k^*)+1}
=\sum_{i=1}^{k^*}s^{2(\mu_{k^*-i+1}-k^*+i-1)+1}
=\sum_{j=1}^{k^*}s^{2(\mu_j-j)+1},
\]
and
\[
\sum_{i=k^*+1}^{k^*+k}s^{-2(w_i+k^*)+1}=\sum_{i=k^*+1}^{k^*+k}
s^{-2(\la_{i-k^*}+k^*-i+1)+1}
=\sum_{j=1}^k s^{-2(\la_j-j)-1}.
\]
Then \[c_{\bf w}=v^{-1}\sum_{j=1}^{k^*}s^{2(\mu_j-j)+1}+v\sum_{j=1}^k
s^{-2(\la_j-j)-1}+\frac{v^{-1}s^{-2k^*}-vs^{2k}}{s-s^{-1}}.\]

Rewrite $s_{\mu,\la}$ by summing over cells in the same row for
each of
$\la$ and $\mu$. Thus
\beqn
s_{\mu,\la}&=&v^{-1}(s-s^{-1})\sum_{x\in\mu}s^{2c(x)}
-v(s-s^{-1})\sum_{x\in\la}s^{-2c(x)}+\delta\\
&=&v^{-1}(s-s^{-1})\sum_{i=1}^{k^*}\sum_{j=1}^{\mu_i} s^{2(j-i)}-
v(s-s^{-1})\sum_{i=1}^k\sum_{j=1}^{\la_i}s^{-2(j-i)}+\de\\
&=&
v^{-1}\sum_{i=1}^{k^*}(s^{2(\mu_i-i)+1}-s^{-2i+1})+v
\sum_{i=1}^{k}(s^{-2(\la_i-i)-1}-s^{2i-1})+\de.
\eeqn

Comparison of the terms in the formulae above then shows that
$c_{\bf w}=s_{\mu,\la}$, noting that 
\beqn
\de-v^{-1}\sum_{i=1}^{k^*}s^{-2i+1}-v\sum_{i=1}^k s^{2i-1}&=&
v^{-1}\frac{s^{-2k^*}}{s-s^{-1}}-v\frac{s^{2k}}{s-s^{-1}}.
\eeqn

We now show that $\qlm\ne 0\in\CC$ by proving
that its framed Homfly evaluation $<\qlm>$ in $\La$ is non-zero.
This calculation uses an explicit formula from Macdonald \cite{Macdonald}
for
$<\qlm>$ in the case where the partition $\mu$ is empty.

There is an expression for $H(t)=\sum_{n=0}^\infty<h_n>t^n$ as an infinite
product
\[H(t)=\prod_{i=0}^\infty \frac{1-bq^it}{1-aq^it}\] with $a=vs,
b=v^{-1}s$ and $q=s^2$,
 \cite{Aiston, AistonMorton2}. For such an
infinite product $H(t)$ Macdonald
\cite{Macdonald} gives the formula
\beqnn
s_\la&=&q^{n(\la)}\prod_{x\in\la}\frac{a-bq^{c(x)}}{1-q^{h(x)}} \label{mac}
\eeqnn  where
$s_\la=<Q_{\la,\phi}>$ is the Schur function for the partition $\la$ given
as a determinant by the Jacobi-Trudy formula from the coefficients
$<h_n>$ of $H(t)$.

 The product in (\ref{mac}) is taken over cells
$x$ of the partition $\la$ with content $c(x)$ and hook length $h(x)$.
Then $<Q_{\la,\phi}>\ne0$ so long as $ab^{-1}\ne q^{c(x)}$ for any cell
$x\in\la$. After the substitution $v=\pm s^N$ this expression is
non-zero in the ring $\La$ if there are no cells $x\in\la$ with $c(x)=N$.
This is the case so long as $N\ge\la_1$.

Given partitions $\la,\mu$ we now choose $N$ and a partition $\nu$ with
$\nu_1\le N$ such that $<\qlm>=<Q_{\nu,\phi}>$ after substituting $v=-s^N$.
It follows that the resulting value of $<\qlm>$ is non-zero, and hence that
$\qlm\ne0\in\CC$.

It is enough to choose $N\ge\la_1+\mu_1$ and take \[\nu_i=\begin{cases}
N-\mu_{k^*-i+1}, & i\le k^*\\
\la_{i-k^*}, &k^* <i\le k^*+k.
\end{cases}\]

In \cite{AistonMorton2} there is an explicit formula 
\[<h_n>=(-1)^n\prod_{i=1}^n \frac{vs^{-i+1}-v^{-1}s^{i-1}}{s^i-s^{-i}}.\]

Then $<h_n>=<h_{N-n}>, 0\le n\le N,$ when $v=-s^N$,  since
\beqn
\prod_{i=1}^N(s^i-s^{-i})&=&
\prod_{i=1}^n(s^i-s^{-i})\prod_{j=1}^{N-n}(s^j-s^{-j})<h_n>\\
&=&\prod_{i=1}^n(s^i-s^{-i})\prod_{j=1}^{N-n}(s^j-s^{-j})<h_{N-n}>
\eeqn for this value of $v$. The result holds also for  $n<0$ and $n>N$.

When we replace every entry $h_n^*$ in the starred rows of a matrix
$M_{\bf w}=\begin{pmatrix}h^*&\dashv\\
h&\dashv\end{pmatrix}_{\bf w}$ by $h_{N-n}$ we get a matrix $M_{\bf
v}=\begin{pmatrix}h&\dashv\end{pmatrix}_{\bf v}$ with standard rows only
and index vector
 given by \beqnn
v_i&=&\begin{cases} N-w_i,&i\le k^*,\\w_i,& i>k^*.
\end{cases}\label{Nsubs}
\eeqnn
  The matrices $<M_{\bf w}>$ and $<M_{\bf v}>$
given by evaluating the framed Homfly polynomial of their entries  are then
identical after the substitution
$v=-s^N$, since $<h^*_n>=<h_n>\in\La$.

When $M_{\bf w}$ is given by  partitions $\la$ and
$\mu$, and $N\ge\la_1+\mu_1$, the matrix $M_{\bf v}$ arises from the
partition $\nu:=\la\cup(N-\mu)$ above with $k^*+k$ parts
$\nu_1=N-\mu_{k^*}\ge N-\mu_{k^*-1}\ge\cdots\ge N-\mu_1\ge\la_1\ge
\cdots\ge\la_k>0$.

Now $<\qlm>=<\det M_{\bf w}>=\det <M_{\bf w}>$ and
$<Q_{\nu,\phi}>=\det<M_{\bf v}>$. Then $<\qlm>=<Q_{\nu,\phi}>$ after substituting
$v=-s^N$, and 
$<Q_{\nu,\phi}>\ne 0$  after this substitution, since $\nu_1\le N$. It
follows that $\qlm\ne0$, and is hence an eigenvector for $\ph$.
\end{proof}

We thus have identified an eigenvector $\qlm$ for each eigenvalue
$s_{\la,\mu}$ of
$\ph$, giving the explicit basis of $\CC$. 

\subsection{Properties of the
 basis elements $\qlm$}
When an element $X\in\CC$ is used to decorate a framed knot $K$, the Homfly
polynomial, $P(K;X)\in\La$, of the decorated knot gives a 2-variable
invariant of
$K$ determined by $X$. Each basis element $\qlm$ gives such an invariant of
$K$, and these determine the invariants for all choices of $X$. When
$\mu=\phi$ these invariants, for each fixed $\la$, give a family of
1-variable invariants after setting
$v=s^{-N}$ which coincide, up to a power of $s$, with the quantum group
invariants of
$K$ determined by the irreducible $sl(N)_q$-module corresponding to the
partition $\la$. 

In general, the 2-variable invariant from $\qlm$, when
evaluated at $v=s^{-N}$, gives the $sl(N)_q$ quantum invariant for an
irreducible module whose partition depends on $\la,\mu$ and $N$. The
simplest example where both $\la$ and $\mu$ occur non-trivially is
$Q_{\Box,\Box}$, where $\Box$ represents the unique partition of $1$. The
2-variable invariant in this case, after setting $v=s^{-N}$, evaluates the
quantum invariant when the knot is coloured by the adjoint representation
of
$sl(N)_q$ for each $N$.    Since
$Q_{\Box,\Box}=h_1h^*_1-1$ this invariant is very nearly the Homfly
polynomial of the reverse parallel of $K$ with one strand in each
direction.

We have already noted that $\overline{\qlm}=\qlm$ and that
$\qlm^*=Q_{\mu,\la}$. Hence the invariant of $K$ when decorated by $\qlm$
is symmetric in $\la$ and $\mu$, and is conjugated in $\La$ when $K$ is
replaced by its mirror image.

In the skein $\CC$ the replacement of $s$ by $-s^{-1}$, while
fixing $v$ and leaving diagrams unaltered, has the effect of interchanging
rows and columns in the Young diagram of each partition $\la$ to give its
conjugate partition $\la^\vee$. Then $\qlm$ becomes
$Q_{\la^\vee,\mu^\vee}$, and the invariant  $P(K,Q_{\la^\vee,\mu^\vee})$
is given  by changing $s$ to $-s^{-1}$ in $P(K:\qlm)$. Where both
partitions
$\la$ and $\mu$ are self-conjugate the invariant is symmetric in $s$ and
$-s^{-1}$, and can then be written in terms of
$v$ and
$z=s-s^{-1}$, as is the case for the ordinary Homfly polynomial of $K$,
which comes from the decoration of $K$ by the single string
$Q_{\Box,\phi}$.

\begin{theorem}
The basis elements $\qlm$ of $\CC$  have the property that the product of
any two is an integer linear combination of basis elements with
non-negative coefficients. \end{theorem}
\begin{proof}
 In the case of the subalgebra $\CC^+$, the basis
restricts to the elements
$Q_{\la,\phi}$ which behave like the basis $\{s_\la\}$ of Schur functions
for the algebra of symmetric functions in many variables. Products of
these are non-negative integer combinations of Schur functions, whose
coefficients can be determined combinatorially by the
Littlewood-Richardson rules. 

The algebra $\CC$ is the polynomial algebra generated by
$\{h_n\},\{h_n^*\}, n\ge1$, with $1=h_0=h^*_0$ and $0=h_n=h^*_n, n<0$.
Define a homomorphism
$\ph_N:\CC\to\CC^+$ by
\[\ph_N(h_n^*)=
h_{N-n},\  0<n, \quad \ph_N(h_n)=\begin{cases}h_n,& 0< n<N,\\ 1,&
n=N,\\0,& n>N.
\end{cases}\] 
 The definition ensures that
$\ph_N(h_n^*)=\ph_N (h_{N-n})$ for {\em all} integers $n$.  Then 
$\ph_N(\qlm)=\ph_N(Q_{\nu,\phi})$ when $N\ge
\la_1+\mu_1$ and $\nu=\la\cup(N-\mu)$ is the partition with parts
$\{\la_i\},\{N-\mu_j\}$.
 The map $\nu\to\ph_N(Q_{\nu,\phi})$ is
injective on partitions $\nu$ with $\nu_1<N$, since the monomial
$h_{\nu_1}h_{\nu_2}\ldots h_{\nu_k}$ can be recovered unambiguously from
the polynomial $\ph_N(Q_{\nu,\phi})$. Indeed the set
$\{\ph_N(Q_{\nu,\phi}):\nu_1<N\}$ forms a basis for $\ph_N(\CC)$.

 For any
finite set of elements
$\qlm$ we can choose $N$ sufficiently large so that $\ph_N$ is injective
on that set. It is enough to take $N>2\max(\la_1,\mu_1)$ for all
$\la,\mu$ in the set, as then the pair $\la,\mu$ can be recovered from
the partition $\la\cup(N-\mu)$ by considering the parts $<N/2$ and those
$>N/2$.

 Expand any product $\qlm Q_{\la',\mu'}$  as a linear combination
$\sum a_{\la'',\mu''}Q_{\la'',\mu''}$ of basis elements of
$\CC$, and choose $N$ so that the supporting basis elements are mapped
 by $\ph_N$ to distinct basis elements of $\ph_N(\CC)$. Write
$\ph_N(Q_{\la'',\mu''})=\ph_N(s_{\nu''})$ with $\nu_1''<N$ for these basis
elements, and also write $\ph_N(Q_{\la,\mu})=\ph_N(s_{\nu})$ and
$\ph_N(Q_{\la',\mu'})=\ph_N(s_{\nu'})$.  
Expand the product $s_\nu s_{\nu'}=\sum b_\rho s_\rho$ with
non-negative integer coefficients $b_\rho$. Then
\[\sum a_{\la'',\mu''}\ph_N(s_{\nu''})=\sum  b_\rho \ph_N(s_\rho).\]
Now $\ph_N(s_\rho)$ is a basis element of $\ph_N(\CC)$ when $\rho_1<N$. It
is also a basis element if $\rho_1=N$, and is zero if $\rho_1>N$.
 A comparison of the two
sides shows that each non-zero coefficient $a_{\la'',\mu''}$ is a sum of
one or more positive integers $b_{\rho}$.
\end{proof}

\begin{remark} This result has a close connection with the behaviour of
irreducible mixed tensor representations of $GL(N)$, which are
characterised by pairs of partitions $\lm$. Work of King
\cite{King}, Stembridge
\cite{Stembridge} and Koike \cite{Koike} 
could be used to give an explicit combinatorial calculation of the
coefficients in a product
$\qlm Q_{\la',\mu'}$.
\end{remark}

\subsection{$(1,1)$-tangle invariants}
For a framed knot $K$ the invariants $P(K:\qlm)$ can be expressed in terms
of the eigenvalues of a map $K:\CC\to\CC$, described as follows. Draw a
diagram of
$K$ as a closed (1,1)-tangle in the annulus and decorate this diagram by
any element
$X\in\CC$ to give the element $K(X)\in\CC$. Thus $P(K:X)$ is the Homfly
polynomial $<K(X)>$. 

Now the map $K$ commutes with the meridian map $\ph$,
and all eigenspaces of $\ph$ are 1-dimensional, so every eigenvector of
$\ph$ is also an eigenvector of $K$. We can then write
$K(\qlm)=K_{\la,\mu}\qlm$ for some scalar $K_{\la,\mu}$. The eigenvalues
of $K$ are 
\[K_{\la,\mu}=\frac{P(K:\qlm)}{<\qlm>},\] giving a normalised version of
$P(K:\qlm)$ which takes the value 1 on the unknot. These eigenvalues are
sometimes known as (1,1)-tangle invariants of $K$, and have the advantage
that they often remain non-zero under  evaluations which send the
corresponding invariant $P$ to zero irrespective of $K$.

When $\mu=\phi$ it can be shown \cite{Morton2}, using properties of the
Gyoja-Aiston idempotents, that the eigenvalues $K_{\la,\phi}$, although
conceivably  rational functions in $v$ and $s$, are in fact integral
elements of $\La$, lying in the subring $\La_0={\bf
Z}[s^{\pm1},v^{\pm1}]$.  This is also the case for $K_{\Box,\Box}$,
and we believe that it is true in general.

By way of example we give the invariant $K_{\Box,\Box}$ for the trefoil
and for the figure-eight knot, in terms of $v$ and $z=s-s^{-1}$.

For the trefoil we get $v^2-4v^4+4v^6
+z^2(1+2v^2-7v^4+4v^6)+z^4(v^2-2v^4+v^6)$ for some choice of framing -
there is a factor of $v^2$ to be used for each change of framing.

For the figure-eight with zero framing the result is symmetric in
$v^{\pm1}$, and is
$3-2z^2-6z^4-2z^6+(v^2+v^{-2})(-2-z^2+2z^4+z^6)+(v^4+v^{-4})(1+2z^2+z^4)$.

Version 1.7, July 2004

Copies of items marked * can be found on the Liverpool Knot Theory site
{\tt http://www.liv.ac.uk/\~{}su14/knotgroup.html} or by following links from
\break
{\tt http://www.liv.ac.uk/maths/} .
\end{document}